# Testing distribution in deconvolution problems

**Denys Pommeret**

*Institut de Mathématiques de Luminy - Case 907 - Université de la Méditerranée - 13288 Marseille cedex 9, France e-mail:* **pommeret@iml.univ-mrs.fr**

**Abstract:** In this paper we consider a random variable $Y$ contaminated by an independent additive noise $Z$. We assume that $Z$ has known distribution. Our purpose is to test the distribution of the unobserved random variable $Y$. We propose a data driven statistic based on a development of the density of $Y + Z$, which is valid in the discrete case and in the continuous case. The test is illustrated in both cases.

## Contents

## 1. Introduction

Consider the convolution

$$X = Y + Z, \qquad (1.1)$$

where $Y$ and $Z$ are independent r.v. with densities $f$ and $g$ w.r.t. a reference measure $\mu$. According to the nature of the variable, $\mu$ can be the Lebesgue measure on $\mathbb{R}$ (or on interval), or the counting measure $\sum \delta_n$, where $\delta_n(x) = 1$ if $x = n$ and 0 otherwise. The error term $Z$ is assumed to have a known density $g$. Observing $X$ instead of $Y$ and $Z$ yields to a deconvolution problem which consists in distinguishing the two components of the variable. Our purpose is to test the distribution of $Y$.

Different authors have considered the nonparametric deconvolution problem of estimating $f$ or its associated distribution function from i.i.d. observations; that is, to recover the distribution of $Y$ using the contaminated measurements $X$. Nonparametric kernel-type estimation of $f$ based on contaminated data was studied by (2), (6) and (3) among others. This problem is also related to mixture problem since (1.1) may be seen as a location mixture $\int f(x-m)\Pi(dm)$, with $m$ the location parameter and $\Pi$ the known mixing distribution (see for instance (4) and (11)).

Recently, (7) proposed a new procedure for testing the density of $Y$. Their method is based on Fourier transforms of variables $Y$, $Z$ and $X$ (via its nonparametric estimator). In this paper we also use a nonparametric representation of the density $h$ of $X$ to test the distribution of $Y$. But our approach is based on a polynomial expansion of the density $h$ under the null hypothesis. The comparison of this density with the empirical one allows to build a test statistic. A data driven approach permits to select automatically the number of components of the statistic. In fact, we consider the same problem than (7), but instead of considering the problem as a $L^2$ distance between a nonparametric deconvolution estimator and a smoothed version of the density $f$, we restrict our attention to a known reference measure $\mu$ and its associated $L^2$ basis where we represent all densities. In (7), the authors need to make different assumptions on Fourier transforms and the use of kernel density estimators requires classical supplementary assumptions on the bandwidth. In the procedure we proposed in this paper we did not need such parameter thanks to the use of a data driven technic. Then our additional assumption is that densities are squared integrable with respect to $\mu$. Also, for asymptotic results we make technical conditions on eigenvalues of the variance matrix. But in practice, the test may be used without data driven procedure for a large enough number of statistics components.

More recently, such data driven approach has been used for deconvolution problem: in (9), the problem of testing the density of $Y$ is considered. The idea is based on the convolution formula obtained in the continuous case; that is, when the reference measure $\mu$ is the Lebesgue one. In that case the author constructed a score test statistic which is combined with a model selection rule. In our paper, the data driven method is close to the work of (9). Also our test statistic can be related to the score one (see Remark 2.1). However we did not use the classical convolution formula, what allows us to consider discrete convolution problems, as that illustrated in our simulations. Also it permits to introduce a dependance between $Y$ and $Z$ as we underline it then. But the extension in a multivariate setting, as done in (7), requires here boring calculus because orthogonal polynomials are not well developed in this frame. Finally, the novelty of our test is that we can also test discrete convolution while in our knowledge it is essentially the continuous case which is studied in the literature. Also, we extend our study to the dependent case, testing the conditional distribution of $Y|Z$.

The paper is organised as follows: in Section 2 we introduce the method based on the polynomial expansion of the (null) density and we propose a data driven procedure for testing the distribution of the contaminated density. In Section 3 a simulation study shows some comportments of the proposed test.

## 2. Method

Let $X_1, \cdots, X_n$ be i.i.d. random variables with density $h$, satisfying the convolution formula (1.1). We assume that $g$ is known and we want to test

$$H_0 : f = f_0 \quad \text{against} \quad f \neq f_0,$$

where $f_0$ is a fixed density. Let $\mu$ be a probability measure on $\mathbb{R}$ with density $m$ and such that both distributions of $Y$ and $X$ are dominated by $\mu$. Denote by $\mathcal{B} = \{Q_i; i = 0, 1, \cdots\}$ an associated basis of dense orthogonal functions with respect to $\mu$. When $h$ belongs to $L^2(\mu)$ (what we shall suppose) we have the following expansion:

$$h(x) \quad = \quad \sum_{n \in \mathbb{N}} a_i Q_i(x), \tag{2.1}$$

where $a_i = \mathbb{E}(Q_i(X)m(X))$. Our method consists in comparing the coefficients $a$'s based on this expansion with those obtained under the null hypothesis. Under $H_0$, equality (2.1) may be rewritten as

$$h(x) \quad = \quad \sum_{i \in \mathbb{N}} \mathbb{E}(Q_i(Y+Z)m(Y+Z))Q_i(x)$$

$$= \quad \sum_{i \in \mathbb{N}} \alpha_i Q_i(x).$$

Testing $H_0$ is equivalent to test

$$H_0 : a_i = \alpha_i \qquad \forall i = 1, 2, \cdots,$$

Under $H_0$ these coefficients $\alpha$'s can be easy to calculate with a good choice of the reference measure and its associated polynomials (see illustrations below). Thus, a natural statistic can be constructed as follows: for some integer $k$, write

$$B_k \quad = \quad (\widehat{b}_1, \cdots, \widehat{b}_k),$$

where

$$\widehat{b}_j \quad = \quad \frac{1}{\sqrt{n}} \Big( \sum_{i=1}^n Q_j(X_i)m(X_i) - \alpha_j \Big).$$

By the Central Limit Theorem, we have the following convergence under $H_0$:

$$U_k = \Sigma_k^{-1/2} B_k^T \quad \longrightarrow_{\mathcal{L}} N(0, I),$$

where $\Sigma_k$ is the $k \times k$ covariance matrix $Var(B_k)$.

Under $H_0$, $T_k = \|U_k\|^2$ is asymptotically Chi-squared distributed with $k$ degrees of freedom and then we reject the null hypothesis for large values of $T_k$. But the arbitrary choice of $k$ is the weakness of this method and we adapt a data driven method to select the number of components in the test statistics. Imitating the work of (10), we consider an increasing sequence of number of components $k(n)$ such that $\lim_{n\to\infty} k(n) = \infty$. Our selection rule is based on the Schwarz criteria. Write

$$S_n = \min\big\{ \underset{1 \leq k \leq k(n)}{\operatorname{argmax}} (T_k - k\log(n)) \big\}. \tag{2.2}$$

Then our data driven test statistic is $T_{S_n}$.

**Proposition 2.1.** *Let $\lambda_{k(n)}$ be the smallest eigenvalue of $\Sigma_{k(n)}$. Assume that $\frac{\log k(n)}{\lambda_{k(n)}} = o_{\mathbb{P}}(\log n)$. Then, under $H_0$, $T_{S_n}$ converges in distribution to a Chi-squared random variable with 1 degree of freedom.*

*Proof.* Under $H_0$, $T_1$ converges to a Chi-squared random variable with one degree of freedom and then we have to show that $\mathbb{P}(S_n = 1)$ tends to zero.

Since $(S_n = k)$ implies $(T_k - k\log(n) \geqslant T_1 - \log(n))$ we have $\mathbb{P}(S_n = k) \leq P(T_k > (k-1)log(n))$. Then we have

$$\mathbb{P}(S_n \geqslant 2) = \sum_{k=2}^{k(n)} \mathbb{P}(S_n = k) \leq \sum_{k=2}^{k(n)} \mathbb{P}\bigg( T_k \geqslant (k-1)\log(n) \bigg). \tag{2.3}$$

Using the fact that $\frac{1}{\lambda_k} = \sup_{X \in \mathbb{R}^{*k}} \frac{X' \Sigma_k^{-1} X}{X' X}$, we obtain

$$T_k = {}^t B_k \Sigma_k^{-1} B_k \leq \frac{\|B_k\|^2}{\lambda_k} \tag{2.4}$$

and

$$\mathbb{P}\bigg( T_k \geqslant (k-1)\log(n) \bigg) \leq \mathbb{P}\bigg( \|B_k\|^2 \geqslant \lambda_k (k-1)\log(n) \bigg). \tag{2.5}$$

As $\Sigma_k$ is an Hilbert Schmidt operator, its trace is bounded by a constant, say $M$, independently of $k$, and we have (under $H_0$)

$$\begin{aligned}
\mathbb{E}(\|B_k\|^2) &= Var(B_k) \\
&= \sum_{i=1}^{k} \big( \mathbb{E}(Q_i(X)m(X))^2 - (\alpha_i)^2 \big) \\
&= Trace(\Sigma_k) \leq M.
\end{aligned}$$

Combining (2.5) with Markov inequality and the above result yields

$$\mathbb{P}\bigg( T_k \geqslant (k-1)\log(n) \bigg) \leq \frac{M}{(k-1)\log(n)\lambda_k}. \tag{2.6}$$

Finally, using a result on harmonic sum we get $\sum_{k=1}^{n} \frac{1}{k} = O_{+\infty}(\log n)$ and we obtain

$$\mathbb{P}(S_n \geqslant 2) \leq \frac{M}{\inf_{1 \leq k \leq k(n)} \lambda_k} \left( \frac{\log k(n)}{\log n} \right).$$

Since matrices $\Sigma_k$ are embedded, $\lambda_k$ is a decreasing sequence and we get the result. □

**Remark 2.1.** *It is clear that our hypotheses can be rewritten $H_0 : \theta = 0$ against $H_1 : \theta \neq 0$, where $\theta = \mathbb{E}(V(k))$ with $V(k) = (Q_i(X)m(X) - \alpha_i)_{i=1,\cdots,k}$. Then $T_k$ coincides with the score statistic if the maximum likelihood estimator $\widehat{\theta}$ of $\theta$ equals the empirical mean of the sample of the $V(k)$'s, that is $\widehat{\theta} = \frac{B_k}{\sqrt{n}}$, as it is the case for instance when the distribution of $V(k)$ belongs to an exponential family.*

## 3. Some particular cases

### 3.1. Continuous positive case

Assume that $Y$ and $Z$ are continuous positive random variables. We choose $\mu$ the exponential distribution with mean 1 and $Q_i = L_{i,1}$ its associated Laguerre orthogonal polynomials (see (1)). In Appendix we recall some basic properties of Laguerre and generalized Laguerre polynomials. Equality (2.1) may be rewritten as

$$
\begin{aligned}
h(x) &= \sum_{i \in \mathbb{N}} \mathbb{E}(Q_i(Y+Z)\exp(-Y)\exp(-Z))Q_i(x) \\
&= \sum_{i \in \mathbb{N}} \big( \sum_{s \leq i} C_{i,s} \mathbb{E}(E_{s,u}(Y)E_{i-s,v}(Z)) \big) Q_i(x),
\end{aligned}
$$

where $E_{s,u}(x) = L_{s,u}(x)\exp(-x)$, $C_{i,s}$ are coefficients given in Appendix and $u, v$ are arbitrar positive reals satisfying $u + v = 1$. Under $H_0$, using the previous decomposition, we have:

$$
\begin{aligned}
\Sigma_{ij} &= \mathbb{E}(Q_i(Y+Z)Q_j(Y+Z)\exp(-2Y)\exp(-2Z)) - \alpha_i\alpha_j \\
&= \sum_{s \leq i}\sum_{t \leq j} C_{s,i}C_{t,j}\mathbb{E}(E_{(s,t),u}(Y)E_{(i-s,j-t),v}(Z)) - \alpha_i\alpha_j,
\end{aligned}
$$

where $E_{(s,t),u}(x) = L_{s,u}(x)L_{t,u}(x)\exp(-2x))$.

### 3.2. Continuous bounded case

Assume that $Y$ and $Z$ are bounded, w.l.o.g. taking values in $[0;1]$. We choose $\mu$ the uniform distribution and $Q_i$ its associated Legendre orthogonal polynomials described in Appendix. Then equality (2.1) may be rewritten as

$$
\begin{aligned}
h(x) &= \sum_{i \in \mathbb{N}} \mathbb{E}(Q_i(Y+Z))Q_i(x) \\
&= \sum_{i \in \mathbb{N}} \big( \sum_{s+t \leq i} C_{i,s,t} \mathbb{E}(Y^s Z^t) \big) Q_i(x),
\end{aligned}
$$

where $C_{i,s,t}$ are coefficients obtained by expansion of $Q_i$. Under $H_0$ :

$$
\begin{aligned}
\Sigma_{ij} &= \mathbb{E}(Q_i(Y+Z)Q_j(Y+Z)) - \alpha_i \alpha_j \\
&= \sum_{s+t \leq i} \sum_{u+v \leq j} C_{i,s,t} C_{j,u,v} \mathbb{E}(Y^{s+u}) \mathbb{E}(Z^{t+v}) - \alpha_i \alpha_j.
\end{aligned}
$$

### 3.3. Discrete case

Assume that $Y$ and $Z$ are discrete random variables taking values in $\mathbb{N}$. We choose $\mu$ the geometric distribution with probability $\mu(x) = p^x(1-p)$ for $x = 0, 1, \cdots, p \in (0,1)$ and $Q_i = M_{i,1}$ its associated Meixner orthogonal polynomials described in Appendix. Equality (2.1) may be rewritten as

$$
\begin{aligned}
h(x) &= \sum_{i \in \mathbb{N}} \mathbb{E}(Q_i(Y+Z)p^Y p^Z (1-p))Q_i(x) \\
&= \sum_{i \in \mathbb{N}} \big( \sum_{s \leq i} C_{i,s} \mathbb{E}(E_{s,u}(Y) E_{i-s,v}(Z)) \big) Q_i(x),
\end{aligned}
$$

where $E_{s,u}(x) = \mathbb{E}(M_{s,u}(x)p^x(1-p)^{1/2})$, $C_{i,s}$ are coefficients given in Appendix and $u, v$ are arbitrar positive reals satisfying $u + v = 1$. Under $H_0$, using the previous decomposition, we have :

$$
\begin{aligned}
\Sigma_{ij} &= \mathbb{E}(Q_i(Y+Z)Q_j(Y+Z)p^{2Y} p^{2Z}(1-p)^2) - \alpha_i \alpha_j \\
&= \sum_{s \leq i} \sum_{t \leq j} C_{s,i} C_{t,j} \mathbb{E}(E_{(s,t),u}(Y) E_{(i-s,j-t),v}(Z)) - \alpha_i \alpha_j,
\end{aligned}
$$

where $E_{(s,t),u}(x) = M_{s,u}(x) M_{t,u}(x) p^{2x}(1-p)$.

### 3.4. Dependent case

When $Y$ and $Z$ are not independent, problem (1.1) can be traited by conditionning with respect to the noise $Z$. In that case, hypotheses concern the conditional distribution of $Y|Z$; that is $H_0 : f_{Y|Z} = f_0(., Z)$. We can use the same approach, replacing coefficients $\alpha_i$ by the quantities

$$
\mathbb{E}(\mathbb{E}(Q_i(Y+Z)m(Y+Z)|Z)).
$$

## 4. Illustrations

In this section we present the results of two simulation studies of our testing procedure. We consider i.i.d. data $X_1, \cdots, X_n$ generated from two convolution models (1.1) satisfying:

**First model** (**Mod1**): $Y$ has exponential distribution with mean 1 and $Z$ is Chi-squared distributed with 1 degree of freedom. Three alternatives are studied:

$\left\{ \begin{array}{l} \textbf{Alt1} : \text{ instead of } (1.1), X \text{ is a mixture with two components :} \\ \quad 50\% \text{ exponential with mean 2 and } 50\% \text{ Chi} - \text{squared with 2 degrees of freedom.} \\ \textbf{Alt2} : \text{ convolution } (1.1) \text{ with both } Y \text{ and } Z \text{ exponential distributed with mean 1.} \\ \textbf{Alt3} : \text{ convolution } (1.1) \text{ with both } Y \text{ and } Z \text{ Chi} - \text{Squared distributed with degree 1.} \end{array} \right.$

**Second model** (**Mod2**): $Y$ has Poisson distribution with mean 1 and $Z$ has Geometric distributed with mean 1. Three alternatives are proposed:

$\left\{ \begin{array}{l} \textbf{Alt4} : \text{ instead of } (1.1), X \text{ is a mixture with two components :} \\ \quad 50\% \text{ Poisson with mean 2 and } 50\% \text{ Geometric with mean 2.} \\ \textbf{Alt5} : \text{ convolution } (1.1) \text{ with both } Y \text{ and } Z \text{ Poisson distributed with mean 1.} \\ \textbf{Alt6} : \text{ convolution } (1.1) \text{ with both } Y \text{ and } Z \text{ Geometric distributed with mean 1.} \end{array} \right.$

It is clear that for Models 1 and 2, the two convolution's components have distributions with relatively close characteristics and we are interested in detecting a confusion between these components (Alternatives 2-3 and Alternatives 5-6) or we are interested in detecting a mixture of these two components instead of a convolution (Alternative 1 and Alternative 4).

For each model and alternative, we compute the test statistic based on a sample size $n = 50, 100$ and $500$ for a theoretical level $\alpha = 5\%$. The empirical level (resp. power) of the test is defined as the percentage of rejection of the null hypothesis over $10000$ replications of the test statistic under $H_0$ (resp. under Alternative). We can see that for Alternative 2, the power is weak for small samples. For Alternative 4 the power is very low. Then the mixture from Alternative 4 and the convolution from Model 2 are very close and it is very difficult to distinguish them.

<div align="center">
Figure 1 here

Figure 2 here
</div>

## 5. Appendix: Orthogonal polynomials

We follow the notation used in (1).

**Laguerre polynomials** Laguerre polynomials $\{L_{n,1}; n \in \mathbb{N}\}$ are defined by their recurrence relation

$$L_{0,1} = 1 \qquad L_{1,1}(x) = 1 - x$$
$$(i+1)L_{i+1,1}(x) = (2i+1-x)L_{i,1}(x) - iL_{i-1,1}(x)$$

They are orthogonal w.r.t. the exponential distribution with density $f(x) = \exp(-x)$ on $\mathbb{R}^+$. For $\alpha > 0$, Generalized Laguerre polynomials $\{L_{n,\alpha}; n \in \mathbb{N}\}$ are orthogonal w.r.t. the Gamma distribution with density $f(x, \alpha) = \exp(-x)x^{\alpha-1}\Gamma(\alpha)^{-1}$. They satisfy

$$L_{0,\alpha} = 1 \qquad L_{1,\alpha}(x) = \alpha - x$$
$$(i+1)L_{i+1,\alpha}(x) = (2i+\alpha-x)L_{i,\alpha}(x) - (i+\alpha-1)L_{i-1,\alpha}(x)$$

Their norms are given by $\|L_{i,\alpha}\|^2 = \alpha^{-2i}(i!)\Gamma(i+\alpha)/\Gamma(\alpha)$. In the exponential case we simply get $\|L_{i,1}\|^2 = (i!)^2$. These polynomials satisfy the following relation that we used in our simulations study (see (5) or (8) for an idea of the proof) for $\alpha = u + v$, $u > 0$, $v > 0$:

$$L_{n,\alpha}(y+z) = \sum_{s \leq n} C_{s,n} u^s L_{s,u}(y) v^{n-s} L_{n-s,v}(z) \tag{5.1}$$

where $C_{s,n} = n!/(s!(n-s)!)$.

**Legendre polynomials** For $\mu$ the uniform distribution on $(0,1)$, associated orthogonal polynomials are (shifted) Legendre polynomials defined by the recurrence relation :

$$P_0 = 1 \qquad P_1(x) = 2x - 1$$
$$(n+1)P_{n+1}(x) = (2n+1)xP_n(x) - nP_{n-1}(x),$$

and sastifying $\|P_n\|^2 = \int_0^1 P_n(x)^2 dx = (2n+1)^{-1}$.

**Meixner polynomials** For $\mu$ the Pascal distribution, $\mu(x) = (1-p)p^x$, $p \in (0,1)$, $x \in \mathbb{N}$, the associated orthogonal polynomials satisfy the following recurrence formula:

$$M_0 = 1 \qquad M_1(x) = 1 - p - x/p$$
$$p/(1-p)M_{n+1}(x) = ((p-1)x + (1+p)n + p)M_n(x) - (1-p)n^2 M_{n-1}(x),$$

with $\|M_n\|^2 = p^{-n}(1-p)^{2n}(n!)^2$. They are particular case of Meixner polynomials $M_{n;b,c}$ with $b = 1$ and $c = p$. These polynomials also satisfy a relation similar to (5.1) that we used for Illustrations:

$$M_{n;b,c}(y+z) = \sum_{s \le n} C_{s,n} M_{s,u,c}(y) M_{n-s,v,c}(z) \qquad (5.2)$$

where $u + v = b$ and $C_{s,n} = n!/(s!(n-s)!)$.

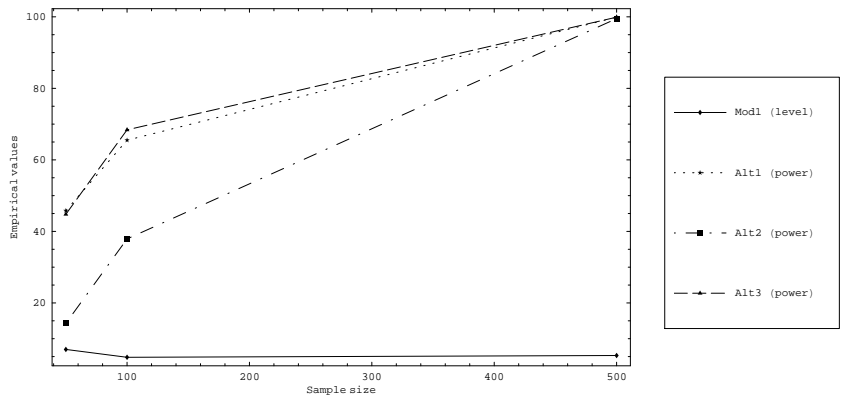

FIG 1. *Empirical level for Model 1 and empirical powers for Alternatives 1-3, with $\alpha = 5\%$, $n = 50, 100, 500$*

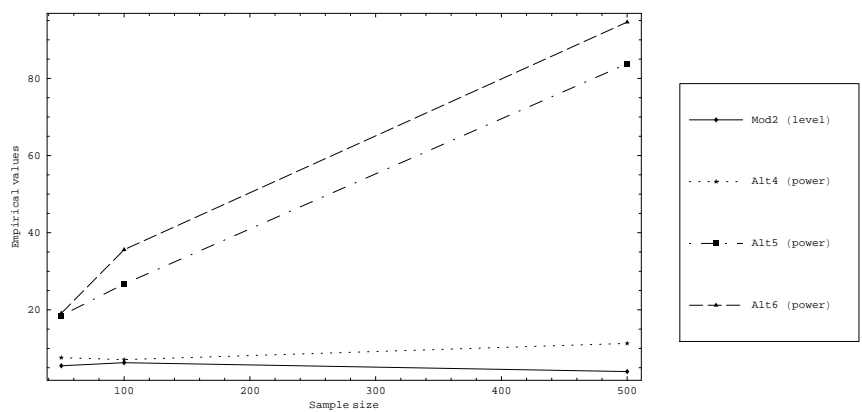

FIG 2. *Empirical level for Model 2 and empirical powers for Alternatives 4-6, with $\alpha = 5\%$, $n = 50, 100, 500$*